\documentclass[12pt]{amsart}

\usepackage{amsmath,amssymb,amscd,amsfonts}

\newtheorem{thm}{Theorem}[section]

\newtheorem{assu-nota}[thm]{Assumption--Notation}
\newenvironment{ipf}{\paragraph{\it Sketch of proof}}{\par\medskip}
\theoremstyle{remark}

\newcommand{\Z}{\mathbb Z}
\newcommand{\Q}{\mathbb Q}

\newcommand{\pp}{\mathbb P}

\newcommand{\Oh}{\mathcal O}

\newcommand{\Si}{\Sigma}

\newcommand{\fie}{\varphi}

\numberwithin{equation}{section}

\begin{document}
\title[The bicanonical map of surfaces  with $p_g=0$ and $K^2\ge 2$]{A survey on the
bicanonical map of surfaces  with $p_g=0$ and $K^2\ge 2$}
\author{Margarida Mendes Lopes \and Rita Pardini}
\date{}

\maketitle

\begin{abstract} We give an up-to-date overview of the known results on the bicanonical map
of surfaces of general type with $p_g=0$  and $K^2\ge 2$.\newline
\noindent 2000 Mathematics Classification: 14J29.
\end{abstract}
\section{Introduction} Many examples of complex surfaces of general type with
$p_g=q=0$ are known, but a detailed classification is still lacking, despite much
progress in the theory of algebraic surfaces. Surfaces of general type are often
studied using properties of their canonical curves. If a surface has $p_g=0$, then
there are of course no such curves, and it is natural to look instead at the
bicanonical system, which is not empty.\par In this survey we describe the present
(December 2001) state of knowledge about the bicanonical map for minimal surfaces
of general type with
$p_g=0$ and $K^2\ge 2$.\par We do not consider the case  $K_S^2=1$ (the so called
{\it numerical Godeaux} surfaces) because, in what concerns the bicanonical map,
this case is special (see \S \ref{dimension}). We just
remark that the numerical Godeaux  surfaces are somewhat  better understood than
the other surfaces of general type with $p_g=0$ and we refer to the paper \cite{cp} and its  bibliography.

This survey is organized as follows: in  section $2$ we discuss the dimension
of the bicanonical image and  in  section
$3$  the base points of $|2K_S|$. In section $4$ we present the bounds on
the degree of the bicanonical map for
$K_S^2\geq 2$. In section $5$ we discuss the  surfaces that occur as bicanonical
images, whilst in section
$6$ we describe  a few relevant examples. Finally in section
$7$ we present some classification results  and in section $8$ we present a list
of  open problems.
\par For each of the results presented we only give a very rough sketch of the
proof, referring to the relevant papers for the missing details.

\bigskip
\noindent{\bf Acknowledgements.}  The present collaboration takes place in the
framework of the european contract EAGER, no. HPRN-CT-2000-00099. The first author
is a member of CMAF and of the Departamento de Matem\'atica da Faculdade de
Ci\^encias da Universidade de Lisboa and the second author is a member of GNSAGA of
CNR and  of the italian  project PIN 2000 ``Geometria delle Variet\`a Algebriche''.
\medskip

\noindent{\bf Notation and conventions.}  We work over the complex numbers; all
varieties are assumed to be  compact and algebraic.  We do not distinguish between
line bundles and divisors on a smooth variety, using the additive and the
multiplicative notation interchangeably. Linear equivalence is denoted by
$\equiv$.   The rest of  the notation is standard in algebraic geometry.

\section{The dimension of the bicanonical image}\label{dimension}
 Let $S$ be a minimal complex surface of general type with $p_g(S)=0$. It is well
known that:
\begin{itemize}
\item $q(S)=0$,
\item $1\le K_S^2\le 9$
\item $P_2(S):=h^0(S, 2K_S)=1+K_S^2$.
\end{itemize}

\medskip
 We denote by $\fie \colon S\to\pp^{K^2_S}$ the bicanonical map of $S$ and by $\Si$
the image of $\fie$. The first question one asks   about the bicanonical map
is what is the dimension of $\Si$. For
$K_S^2=1$, one has
$P_2(S)=2$ and so $\Si$ is a curve. For  $K_S^2\geq 2$, the
answer  was given  by Xiao Gang in the mid-eighties:
\begin{thm} {\em (Xiao Gang, \cite{xiaocan})}\label{Xiaoimage} Let $S$ be a minimal
complex surface of general type with $p_g(S)=0$.  If
$K^2_S\ge 2$ then the image of the bicanonical map of $S$
 is a surface.
\end{thm}

\begin{ipf} We just explain the main ideas and refer the reader to   \cite{xiaocan}
for details.

By contradiction, suppose that $|2K_S|$ is composed with a pencil. Then one can
show that  necessarily
$2K_S\equiv aF+Z $, where $|F|$ is a base point free genus $2$ pencil  and
$a=K_S^2$. Using the very precise description of  Horikawa (\cite{ho}) for the
reducible fibres of a genus $2$ fibration in terms of
$K_S^2$ and
$\chi(\Oh_S)$,  one shows that for   $K_S^2\geq 3$ the surface   $S$ does not
contain a genus $2$ fibration, since otherwise  the components of the reducible fibres
would give  $b_2(S)$ or more independent classes  in
$H^2(S,\Z)$. For $K_S^2=2$, $S$ can have a genus  $2$ fibration with general fibre
$F$,  but one can use   the same type of argument to show that  it is impossible to
decompose $2K_S$ as $2K_S\equiv 2F+Z$, i.e.
 that
$|2K_S|$ is not composed with $|F|$.
\qed\end{ipf}
\medskip So the surfaces with $K_S^2=1$ and $p_g=0$, the  {\it numerical Godeaux},
are in a class of their own.  We  just mention here that there is intensive work in
progress on this subject  by F. Catanese and R. Pignatelli and by  Y. Lee, using in
particular  the bicanonical fibration. As already mentioned in the introduction, for more
facts on numerical Godeaux one can see the paper
\cite{cp}, which has also a very complete list of references.

\section{The base points of the bicanonical system}\label{secbpoints} Recall that,  while
for
$p_g(S)>0$ the bicanonical map is defined at every point of $S$ (\cite{bo}, \cite{red}, \cite{francia}, \cite{cc}, cf.
\cite{ciro}), for
$p_g(S)=0$ it  is still unknown whether $\fie$ is always a morphism. For
$K_S^2\geq 5$ we have:
\begin{thm} {\em (Reider, \cite{red})} Let  $S$ be a minimal surface of general type
with  $p_g=0$ and let
$\fie\colon S\to\pp^{K^2_S}$ be the bicanonical map.

If $K_S^2\geq 5$, then $\fie$ is a morphism.
\end{thm}

\medskip
\noindent{\sl Remark 1\/}  This  is a particular case of Reider's theorem
(\cite{red})  about adjoint systems, which only applies if $K_S^2\geq 5$.

\medskip
\noindent{\sl Remark 2\/} As far as we know, for all the known examples of surfaces
of general type with  $2\leq K_S^2\leq 4$ and $p_g=0$ the bicanonical map is a
morphism.

\medskip For $K_S^2=4$, Lin Weng (\cite{linweng}) has proven that  the base locus
of the bicanonical system contains no $-2-$curve.  This result has later been
improved by Langer:
\begin{thm} {\em(Langer, \cite{langer})} Let $S$ be a minimal surface
 of general type with $K^2_S=4$ and $p_g(S)=0$. Then the system $|2K_S|$ has no
fixed component.
\end{thm}

Still in the case $K^2_S=4$, F. Catanese  and F. Tovena (\cite{ct}) and D. Kotschick
(\cite{kotschick}) have related the existence of base points of the bicanonical system to
properties of the fundamental group of the surface. Since the statements are quite
technical, we just quote here the following consequence of their results:
\begin{thm} {\em (Catanese-Tovena, \cite{ct}, Kotschick,  \cite{kotschick})}
 Let $S$ be a minimal surface
 of general type with $K^2_S=4$ and $p_g(S)=0$.  If $H^2(\pi_1(S), \Z_2)=0$, then
the  bicanonical system
$|2K_S|$ is base point free.
\end{thm}

\section{The degree}

Once one knows that for $K^2\ge 2$ the bicanonical image of a surface $S$ of
general type with $p_g=0$  is a surface, it is natural to look for bounds on the
degree
$d$ of the bicanonical map $\fie$. \par If $K_S^2=2$, the bicanonical image is
$\pp^2$ and therefore $d\geq 2$. On the other hand, $(2K_S)^2=8$ implies
$\deg\fie\leq 8$, , with equality holding if and only if
$\fie$ is a morphism. All the known examples with $K_S^2= 2$ have
$\deg\fie=8$.

 For higher values of $K_S^2$ we have:

\begin{thm}\label{degree1} {\em (\cite{marg})} Let $S$ be a minimal complex surface of
general type such that $p_g(S)=0$, $K_S^2\geq 3$ and let
$\fie\colon S\to
\pp^{K^2_S}$ be the bicano\-ni\-cal map of $S$. Then the degree of $\fie$ is at most
$5$.\par

If $\fie$ is a morphism (in particular if $K_S^2\geq 5$)  then the degree of $\fie$
is at most
$4$.
\end{thm}
\begin{ipf} (See \cite{marg} for the complete proof).
Let $d$ be the degree of $\fie$ and let $m$ be the degree of the bicanonical
image $\Si\subset \pp^{K_S^2}$.  Since $\Si$ is a non-degenerate surface in
$\pp^{K_S^2}$, one has
$m\geq K_S^2-1$.
Write $|2K_S|=|M|+F$, where $M$ and $F$
are the moving part and the fixed part of the system, respectively. Notice
that, if
$F\ne 0$,   then $M^2< (2K_S)^2$, by the 2-connectedness of the
bicanonical divisors. So we have $4K^2_S\ge md$ and equality holds if and only if
$\fie$ is a morphism.
 By an easy calculation we see that to prove  the theorem it is
enough to exclude the possibilities
$K_S^2=5$, $d=5$,  and $K_S^2=3$, $d=6$.\par This  is done by using the
classification of surfaces of degree
 $n-1$ in $\pp^n$ (see \cite{nagata}) to find the possibilities for
$\Sigma$. Then,  using the geometry of $\Sigma$, one is able to build irregular double covers of $S$, which in turn, by a
theorem of De Franchis (\cite{defra}, see also \cite{cetraro}), yield special fibrations on $S$.   Finally, with different
``twists'' for each case, the existence of such a fibration  leads to a contradiction.
\qed \end{ipf}

\medskip
\noindent{\sl Remark\/}
  As mentioned above, for all known   examples of surfaces with $K_S^2>1$
the bicanonical map  is a morphism, and so  the bound $5$ of the theorem  above may
not be effective for
$K_S^2=3,4$. \par  On the other hand, the bound $4$, if $\fie$ is a morphism,  is
effective, as shown by the Burniat surfaces with
$K_S^2=3,...,6$ (\cite{burniat}, \cite {peters}, see
also Example 3 of \S \ref{examples}).

\bigskip For high values of $K_S^2$ these bounds can be improved.
\begin{thm}\label{degree2} {\em (\cite{mp})} Let $S$ be a minimal complex surface of
general type such that $p_g(S)=0$ and let
$\fie\colon S\to
\pp^{K^2_S}$ be the bicanonical map of $S$. Then one has the following bounds on
$d:=\deg\fie$:
\begin{itemize}
\item[i)] if $K_S^2=9$, then $d=1$;
\item[ii)] if $K_S^2=7,8$, then $d\le 2$;

\end{itemize}
\end{thm}
\begin{ipf} (See \cite{mp} for the proof). For  $K_S^2=9$, one has
$c_2(S)=3$ and so
$b_2(S)=1$.  The assertion is proven  by combining  Reider's theorem (\cite{red}) and
the fact that
$H^2(S,\Q)$  is  generated by the class of an ample divisor $D$ with $D^2=1$.

For  $K_S^2=7$,  $\fie$ is a morphism, and therefore  the degree $d$ of
$\fie$ is  either
$1$,
$2$, or
$4$. If $d=4$, then  the bicanonical image $\Si$ is a linearly normal surface  of degree $7$
in
$\pp^7$ with
$p_g=q=0$, and so it is the anti-canonical image of
$\pp^2$ blown-up at two points $P,Q$.  Combining the information obtained from the
geometry  of $\Si$ and the fact that the second Betti number of $S$  is small ($b_2(S)=3$),
it is possible to find a contradiction, which shows that $d=4$ does not occur.

For $K_S^2=8$, the technique of proof is analogous.
\qed \end{ipf}

\noindent{\sl Remark\/}  The bounds in this theorem are effective, since there are
examples of minimal surfaces with $p_g=0$ and
$K_S^2=7,8$ for which the bicanonical map has degree $2$ (see Examples $1$ and $2$
in Section \ref{examples}).

\section{The image}

Another natural question that arises is  finding the possibilities for  the image of the bicanonical map $\fie$,
if $\fie$ is not birational.  A priori, one only knows that the image of $\fie$ is a surface with $p_g=q=0$. It turns out that
it is possible to be more precise, as we will see in the next theorem.
\begin{thm}\label{image} {\em(Xiao Gang, \cite{xiao2}, \cite{mp3})} Let $S$ be a minimal complex
surface of general type such that $p_g(S)=0$ and $K_S^2\ge 2$ and let
$\fie\colon S\to \Si\subset\pp^{K^2_S}$ be the bicanonical map of $S$. If $\fie$ is
not birational, then either

i) $\Si$ is a rational surface,

 or

ii) $K^2_S=3$,
$\fie$ is a morphism of degree 2 and $\Si\subset \pp^3$ is an Enriques sextic.
\end{thm}

 \begin{ipf} (See \cite{xiao2} and \cite{mp3} for the proof). If the degree $d$ of $\fie$ is bigger than $2$,
then
$\Si$ is a linearly normal  surface with
$p_g=q=0$ and degree lesser than  or equal to
$2n-2$ in $\pp^n$, and so  a rational surface. In \cite{xiao2}, Xiao Gang, using
double covers techniques, showed that if
$d=2$ and $\Si$ is not rational, then $K_S^2\leq 4$ and $\Si$ is birationally
equivalent to an Enriques surface with $K_S^2+4$ nodes. The  Enriques surfaces with 8 nodes
are classified in \cite{mp3}. Using the knowledge of the linear systems on these surfaces, one
is able to determine precisely the  surfaces of general type with $p_g=0$ and
$K^2=4$ whose bicanonical map factors through a degree 2 map onto an Enriques surface. Such
surfaces had been previously constructed by D. Naie (\cite{naie}) and their bicanonical map
is of degree 4 onto a rational surface.
Hence  the possibility   $K_S^2=4$ is excluded.
\qed\end{ipf}

\section{Some examples}\label{examples}

In this section we give a quick description of some of the known examples of
surfaces with
$p_g=0$.
\bigskip

\noindent{\bf Example 1:} {\em Surfaces with $K_S^2=8$.} All the examples known to
us of  surfaces $S$ with $p_g=0$ and $K_S^2=8$ are obtained by the following
construction, first suggested by Beauville (cf.  \cite{bv}, \cite{dolg}). One takes curves
$C_1$,
$C_2$ of genera $g_1$ and $g_2$, respectively, such that there exists a group $G$
of order
$(g_1-1)(g_2-1)$ that acts faithfully on both $C_1$ and $C_2$. If the quotient
curves
$C_1/G$ and $C_2/G$ are rational and the diagonal action of $G$ on $C_1\times C_2$
is free, then $S:=(C_1\times C_2)/G$ is a minimal surface of general type with
$p_g=0$ and $K_S^2=8$. By \cite{pianidoppi}, the  surfaces with these invariants
and     bicanonical map of degree 2 are obtained by   this construction taking
$C_1$  hyperelliptic of genus 3, and they belong to four different families.
Examples with birational bicanonical map do exist. For instance,  if one  takes
$C_1=C_2$ to be the Fermat quintic, then  it is possible to let $G=\Z_5^2$ act on
the two copies of the curve in such a way that the  diagonal action is free. The resulting surface has birational
bicanonical map.
\medskip

\noindent{\bf Example 2:} {\em A surface with $K_S^2=7$.} This example is due to
Inoue (\cite{inoue}), who costructs it by taking the quotient of a complete
intersection inside the product of four elliptic curves  by a group isomorphic to
$\Z_2^5$ acting freely. Alternatively, this surface can be constructed as a
$\Z_2^2-$cover of a rational surface with 4 nodes (see \cite{mp}).  The bicanonical
map has degree 2 and the bicanonical involution belongs to the Galois group of the
cover. The quotient of this surface by the bicanonical involution is a rational
surface with 11 nodes.
\medskip
\begin{figure}[ht]
\setlength{\unitlength}{0.5cm}
\centerline{\begin{picture}(24,14)
\thicklines
\put(0,0){\line(1,0){24}}
\put(0,0){\line(1,1){12}}
\put(24,0){\line(-1,1){12}}
\thinlines
\put(0,0) {\line(2,1){14}}
\put(0,0){\line(3,1){14}}
\put(24,0){\line(-3,1){14}}
\put(24,0){\line(-2,1){14}}
\put(12,12) {\line(-1,-4){2.4}}
\put(12,12) {\line(1,-6){1.7}}
\put(12.5,12.5){$P_3$}
\put(0,0.9){$P_1$}
\put(24,0.9){$P_2$}
\put(10.5,10){$m_1^3$}
\put(12.4,10){$m_2^3$}
\put(2.2,0.2){$m_1^1$}
\put(2.6,2.2){$m_2^1$}
\put(4.5,5.75){$l_{13}$}
\put(12,0.3){$l_{12}$}
\put(19,5.75){$l_{23}$}
\put(19,0.4){$m_1^2$}
\put(19,2.7){$m_2^2$}
\end{picture}}
 \caption{The branch locus of the Burniat surfaces in $\pp^2$} \label{fig!quad1}
 \end{figure}

\begin{figure}[ht]
\setlength{\unitlength}{0.5cm}
\centerline{\begin{picture}(24,14)
\thicklines
\put(6,0){\line(1,0){12}}
\put(2,5){\line(1,1){6}}
\put(22,5){\line(-1,1){6}}
\thinlines
\put(6,10) {\line(1,0){12}}
\put(7,0) {\line(-1,2){4}}
\put(17,0) {\line(1,2){4}}
\put(11,12){\line(0,-1){6}}
\put(13,12){\line(0,-1){6}}
\put(4,4){\line(1,1){4}}
\put(4,2){\line(1,1){4}}
\put(20,4){\line(-1,1){4}}
\put(20,2) {\line(-1,1){4}}
\put(9.5,10.5){$e_3$}
\put(6.3,1.8){$e_1$}
\put(17,1.8){$e_2$}
\put(10,9){$m_1^3$}
\put(12,9){$m_2^3$}
\put(5,6){$m_1^1$}
\put(4.9,4.2){$m_2^1$}
\put(4.5,8.5){$e'_2$}
\put(12,0.3){$e'_3$}
\put(19,8.5){$e'_1$}
\put(18,4.2){$m_1^2$}
\put(18,6){$m_2^2$}
\end{picture}}
 \caption{The branch locus of the Burniat surfaces in  $\Sigma$ } \label{fig!quad}
 \end{figure}
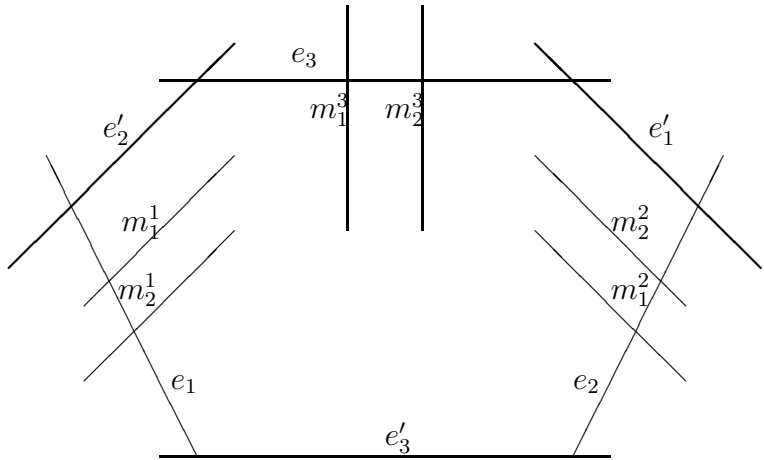

\noindent{\bf Example 3:} {\em Burniat surfaces.} These surfaces were discovered  by
Burniat (\cite{burniat}) and studied later by Peters (\cite{peters}). They are
obtained by taking
$\Z_2^2-$covers of the plane branched on a configuration of lines as shown in Figure
\ref{fig!quad1} in such a way that the images of the divisorial components of the fixed
loci of the three nonzero elements of the Galois group of the cover are $l_{12}+m^1_1+m^1_2$,
$l_{23}+m^2_1+m^2_2$ and
$l_{13}+m^3_1+m^3_2$. For a general choice of the  lines $m^i_j$,  the resulting
surface
$Y$ is singular above the points $P_1,P_2,P_3$ and the minimal resolution of $Y$ is
obtained by  taking base change with the blow--up $\hat{\pp}\to{\pp}^2$ of the plane
at
$P_1,P_2,P_3$ and then normalizing. In this way one obtains a
minimal surface $S$ with
$p_g=0$ and
$K_S^2=6$. The bicanonical map is the composition of the induced
$\Z_2^2-$cover $S\to \hat{\pp}$ with the embedding of  $\hat{\pp}$
as a Del Pezzo sextic in $\pp^6$.

Examples with the same properties and with $2\le K_S^2\le 5$ can be obtained by
letting one or more subsets of 3 lines $m^i_j$  go through the same point.
\section{The limit cases}

Surprisingly, the limit cases of Theorem \ref{degree1} and  Theorem \ref{degree2}, namely
$K^2_S=7,8$, $\deg\fie=2$, and $K^2_S=6$,
$\deg\fie=4$, can be described precisely, as we will
see in the next two theorems.

\begin{thm}\label{k78} {\em (\cite{mp4})} Let $S$ be a minimal complex surface of general
type such that $p_g(S)=0$ and $K^2_S=7,8$.   Let
$\fie\colon S\to
\pp^{K^2_S}$ the bicanonical map of $S$ and assume that $\deg\fie=2$.  Then $K_S$
is ample and:
\begin{itemize}

\item[i)] there exists a fibration
$f\colon S
\to\pp^1$ such that the general fibre $F$  of $f$ is a smooth hyperelliptic curve
of genus $3$ and the bicanonical involution induces the hyperelliptic involution on
$F$;

\item[ii)] if   $K^2_S=8$ then  $f$ is isotrivial and the singular fibres of $f$
are $6$ double fibres with smooth support, while if $K^2_S=7$ then $S$ has $5$
double fibres and exactly one fibre with reducible support.

\end{itemize}
\end{thm}
\vskip0.2truecm\noindent{\sl Remark 1\/}  For $K_S^2=8$, the fact that $f$ is an isotrivial
fibration whose only singular fibres are double fibres with smooth support
implies that $S$ is one of the Beauville surfaces (see \S \ref{examples}, Example 1). Using
this fact it is possible to give a complete classification of these surfaces (see
\cite{pianidoppi}). They belong to 4 different types, and the surfaces of each type form an
irreducible connected component of the moduli space of surfaces of general type.
An interesting feature of these surfaces is that they are smooth minimal models of
 double covers of the plane branched on a curve with certain singularities, a
construction that  had been suggested by Du Val (\cite{duval}, see also \cite{ciro}). The expected number of
parameters of the branch curve of this construction is negative, hence it seems
very difficult prove directly its existence, that  instead  follows ``a posteriori''
from the classification of \cite{pianidoppi}.
\vskip0.2truecm\noindent{\sl Remark 2\/} For $K_S^2=7$ the hyperelliptic fibration is not
isotrivial and a complete classification seems out of reach.  In the Inoue's surface (see
\S \ref{examples}), which is the only known example with $K_S^2=7$ and $\deg\fie=2$, the
unique  fibre with reducible support of the fibration
$f$ is one of the double fibres. In principle one would expect  this to  be a special
situation, hence it would be interesting to find examples where the reducible fibre is not a
double fibre.
\vskip0.3truecm \begin{ipf} (See \cite{mp4} for the proof). Consider the quotient
$Y$ of  $S$ by the bicanonical involution
$\sigma$. $Y$ is a rational surface whose only singularities are
 $\nu=K_S^2+4$ nodes,  which correspond to the isolated fixed points of $\sigma$.
The minimal resolution of
$Y$ is a rational surface $X$ having $\nu\geq b_2(X)-3$  disjoint  $-2-$curves.
Such surfaces are characterized in
\cite{nodi}, where it is shown in particular that there exists a fibration $g\colon Y\to\Si$
with rational fibres and with $[\frac{\nu}{2}]$ double fibres. Now, using  some geometrical
reasoning,  one shows that $g$ pulls back to a fibration $f\colon S\to\pp^1$ such that the
general fibre of $f$ is hyperelliptic of genus 3.
\par The fact that
$K_S$ is ample follows trivially for $K_S^2=8$ from Miyaoka's results (\cite{miyaoka}) on the
existence of rational curves on surfaces. For $K_S^2=7$, the non-existence of $-2-$curves on
$S$ is obtained by analyzing the structure of the unique reducible fibre of $f$ and using the
equality
$b_2(S)=3$.
\qed\end{ipf}
\begin{thm}\label{burniat} {\em(\cite{mp2})} Let $S$ be a minimal complex surface of general
type such that $p_g(S)=0$ and $K^2_S=6$  and let
$\fie\colon S\to
\pp^{K^2_S}$ the bicanonical map of $S$. Then:

\begin{center} $\deg\fie=4$ if and only if $S$ is a Burniat surface.\end{center}

In particular, $K_S$ is ample.
\end{thm}
\begin{ipf} (See \cite{mp2} for the proof).
 The first step in the proof of this theorem consists in showing that  the
bicanonical image
$\Sigma$ of
$S$ is the non-singular Del Pezzo surface of degree $6$ in $\pp^6$. This is shown
by a case by case  exclusion  of all the possible singular such
 $\Sigma$.\par The second step consists in showing that the
sides of the   ``hexagon'' of $-1-$curves   of $\Sigma$ (cf. Fig. 2) are in
the branch locus of
$\fie$.   Using the curves in $S$ which correspond to these sides, one produces a
subgroup
$H$ of Pic$(S)$ such that $H\simeq
\Z_2^3$.
\par By studying the \'etale covers of $S$ given by the nonzero elements of $H$, it is
possible to show that the three pencils
 of $\Si$ corresponding to the lines through $P_1,P_2$ and $P_3$ pull-back in $S$ to
three genus
$3$ hyperelliptic pencils, each having   two irreducible  double fibres beside the two
reducible ones given by pairs of sides of the hexagon. The images of these
irreducible fibres are the remaining components of the branch locus of $\fie$ (see Fig. 2).
\par
Finally, one verifies that  the bicanonical map is composed with the  three involutions of $S$
induced by the hyperelliptic pencils.  It follows that the  bicanonical  map is a
$\Z_2\times\Z_2-$cover and
$S$ is a Burniat surface.
\qed\end{ipf}

Using Theorem \ref{burniat},  one can  obtain very precise information on the geometry of
the moduli space of surfaces with $p_g=0$, $K_S^2=6$ and bicanonical map of degree 4.
 \begin{thm}\label{main2} {\em (\cite{mp2})} Smooth minimal surfaces of general type $S$ with
$K^2_S=6$, $p_g(S)=0$ and bicanonical map of degree $4$ form an unirational
$4-$dimensional irreducible connected component of the moduli space of
surfaces of general type.
\end{thm}
\begin{ipf} (See \cite{mp2} for the proof). By the semicontinuity of $\deg\fie$ and by Theorem \ref{degree1}, the surfaces
with $\deg\fie=4$ are a closed subset of the moduli space.

 Using the theory of natural
deformations of abelian covers, one
constructs explicitly a smooth family  ${\mathcal X}\to B$ of  $\Z_2^2-$covers of the plane
blown up at three non collinear points with the following properties:
\begin{itemize}
\item[i)] $B$ is smooth and  irreducible;

\item[ii)] for every $b\in B$ the fibre $X_b$ is a Burniat surface and every Burniat surface
occurs as a fibre for some $b\in B$;

\item[iii)] the family ${\mathcal X}$ is complete at every point $b$ of $B$.
\end{itemize}

This shows  that the Burniat surfaces are an irreducible  open subset of the moduli space. Hence,
in view of Theorem \ref{burniat}, the surfaces with $\deg\fie=4$ are an irreducible  open and
closed  subset of the moduli space.
\qed\end{ipf}
\noindent{\it Remark 3\/} The limit cases for the degree of the bicanonical map have some
common properties. First of all, by Theorem \ref{k78} and Theorem \ref{burniat}, they all have
ample canonical class. In addition, all surfaces with $K_S^2=6$ and $\deg\fie=4$, all surfaces
with $K_S^2=8$ and $\deg\fie=2$ (see \cite{pianidoppi}) and the Inoue surfaces with $K_S^2=7$
  move in
positive dimensional families, while the expected dimension of the moduli at the
corresponding points is zero.
\section {Some questions}
Here we point out  some  questions that arise naturally from the results outlined in the
previous sections.
\medskip

\noindent{\it Question 1\/} (cf. \S \ref{secbpoints}) Is the bicanonical map
$\fie$  of surfaces with $p_g=0$ a morphism also for $2\le K_S^2\le 4$?
\medskip

\noindent{\it Question 2\/} (cf. Theorem \ref{degree1}) Is there a surface with $p_g=0$,
$K_S^2=3$ or $K_S^2=4$ and $\deg \fie=5$? Notice that for such a surface $\fie$
cannot be a morphism.
\medskip

\noindent{\it Question 3\/} (cf. Theorem \ref{degree1} and Theorem \ref{degree2})
Is there is a surface with $p_g=0$,
$K_S^2=6$ and $\deg \fie=3$?
\medskip

\noindent{\it Question 4\/} (cf. Theorem \ref{burniat}) Is it possible to characterize
surfaces with
$K_S^2=5$,
$p_g=0$  and
$\deg
\fie=4$?
\medskip

\noindent{\it Question 5\/} (cf. \S 7, Remark 2) Are the Inoue surfaces the only surfaces
with $p_g=0$, $K^2_S=7$ and non birational bicanonical map?

\bigskip

\begin{tabbing} 1749-016 Lisboa, PORTUGALxxxxxxxxx\= 56127 Pisa, ITALY \kill
Margarida Mendes Lopes \> Rita Pardini\\ CMAF \> Dipartimento di Matematica\\
 Universidade de Lisboa \> Universit\a`a di Pisa \\ Av. Prof. Gama Pinto, 2 \> Via
Buonarroti 2\\ 1649-003 Lisboa, PORTUGAL \> 56127 Pisa, ITALY\\
mmlopes@lmc.fc.ul.pt \> pardini@dm.unipi.it
\end{tabbing}

\end{document}